%------------------------------------------------------------------------------
% rx.tex  
%                    
% Beginn: 10.11.2011
%------------------------------------------------------------------------------
\magnification=\magstep1   
\input amstex
\UseAMSsymbols
\input pictex
%\hoffset=0truecm \voffset=0truecm 
\vsize=23truecm
\NoBlackBoxes
\parindent=18pt
  
   \font\rmk=cmr8    \font\itk=cmti8  \font\ttk=cmtt8

%\hrule height 2pt \vskip 3pt \hrule \bigskip\bigskip\bigskip

\def\Hom{\operatorname{Hom}}

\def\Ext{\operatorname{Ext}}

\def\rad{\operatorname{rad}}

\def\soc{\operatorname{soc}}
\def\Tr{\operatorname{Tr}}

\def\arr#1#2{\arrow <1.5mm> [0.25,0.75] from #1 to #2}

	\vglue1truecm
\plainfootnote{}
{\rmk 2010 \itk Mathematics Subject Classification. \rmk 
Primary 
        16D90, %Module categories [See also 16Gxx, 16S90]; module theory in a 
               % category-theoretic context; Morita equivalence and duality
        16G10. %Representations of artinian rings 
Secondary:
        16G70. % Auslander-Reiten sequences (almost split sequences) and Auslander-Reiten quivers
}

\centerline{\bf On radical square zero rings.}
                     \bigskip
\centerline{Claus Michael Ringel and Bao-Lin Xiong}     
		  \bigskip\medskip

{\narrower\narrower
Let $\Lambda$ be a connected left artinian ring with radical square zero and with
$n$ simple modules. If $\Lambda$  is not self-injective, 
then we show that any module $M$ with $\Ext^i(M,\Lambda) = 0$ 
for $1 \le i \le n+1$ is projective. We also determine the structure of 
the artin algebras with radical square zero and $n$ simple modules
which have a non-projective module $M$
such that $\Ext^i(M,\Lambda) = 0$ for $1 \le i \le n$.
\par} 
	\bigskip
Xiao-Wu Chen [C] has recently shown: given a connected artin algebra
$\Lambda$ with radical square zero then either $\Lambda$ is self-injective or else 
any CM module is projective.
Here we extend this result by showing: If $\Lambda$ is 
a connected artin algebra with radical square zero and $n$ simple modules
then either $\Lambda$ is self-injective or else any module $M$ with $\Ext^i(M,\Lambda) = 0$ 
for $1 \le i \le n+1$ is projective. 
Actually, we will not need the assumption on $\Lambda$ to be an artin algebra; it is
sufficient to assume that $\Lambda$ is a left artinian ring.
And we show that for artin algebras 
the bound $n+1$ is optimal by
determining the structure of those artin algebras 
with radical square zero and $n$ simple modules
which have a non-projective module $M$
such that $\Ext^i(M,\Lambda) = 0$ for $1 \le i \le n$.
	\medskip
 From now on, let $\Lambda$ be a left artinian ring with radical square zero, this means that
$\Lambda$ has an ideal $I$ with $I^2 = 0$ (the radical) such that $\Lambda/I$ is semisimple
artinian. We also assume that $\Lambda$ is connected (the only central idempotents are $0$
and $1$). The modules to be considered are usually finitely generated left $\Lambda$-modules. 
Let $n$ be the number of (isomorphism classes of) simple modules.

Given a module $M$, we denote by $P M$ a projective cover, by $Q M$ an injective envelope
of $M$. Also, we denote by $\Omega M$ a syzygy module for $M$, this is the kernel
of a projective cover $PM \to M.$ Since $\Lambda$ is a ring with radical square zero,
all the syzygy modules are semisimple. Inductively, we define $\Omega_0 M = M,$ and
$\Omega_{i+1}M = \Omega(\Omega_{i}M)$ for $i\ge 0.$
	\medskip
{\bf Lemma 1.} {\it If $M$ is a non-projective module with $\Ext^i(M,\Lambda) = 0$ for 
$1 \le i \le d+1$ (and $d\ge 1$), then there exists a simple non-projective module
$S$ with $\Ext^i(S,\Lambda) = 0$ for $1 \le i \le d.$}
	\medskip
Proof: We have $\Ext^i(M,\Lambda) \simeq \Ext^{i-1}(\Omega M,\Lambda)$, for all $i \ge 2.$
Since $M$ is not projective, $\Omega M\neq 0$. Now $\Omega M$ is semisimple. If all simple
direct summands of $\Omega M$ are projective, then also $\Omega M$ is projective, but then 
the condition $\Ext^1(M,\Lambda) = 0$ implies that
$\Ext^1(M,\Omega M) = 0$ in contrast to the existence of the exact sequence $0 \to \Omega M
\to PM \to M \to 0.$ Thus, let $S$ be a non-projective simple direct summand of $\Omega M.$
	\medskip
{\bf Lemma 2.} {\it If $S$ is a non-projective simple module with $\Ext^1(S,\Lambda) = 0$,
then $PS$ is injective and $\Omega S$ is simple and not projective.}
	\medskip
Proof: First, we show that $PS$ has length 2. Otherwise, $\Omega S$ is of length at least $2$,
thus there is a proper decomposition $\Omega S = U\oplus U'$ and then there is a canonical
exact sequence
$$
  0 \to PS \to PS/U \oplus PS/U' \to S \to 0,
$$
which of course does not split. But since $\Ext^1(S,\Lambda) = 0$, we have $\Ext^1(S,P) = 0$,
for any projective module $P$. Thus, we obtain a contradiction.

This shows also that $\Omega S$ is simple. Of course, $\Omega S$ cannot be projective,
again according to the assumption that $\Ext^1(S,P) = 0$,
for any projective module $P$. 

Now let us consider the injective envelope $Q$ of $\Omega S$. It contains $PS$ as a
submodule (since $PS$ has $\Omega S$ as socle). Assume that $Q$ is of length at
least $3$. Take a submodule $I$ of $Q$ of length 2 which is different from $PS$ and let
$V = PS+I$, this is a submodule of $Q$ of length 3. Thus, there are the following 
inclusion maps $u_1,u_2,v_1,v_2$:
$$
\CD
\Omega S @>u_1>> PS \cr
     @Vv_1VV          @VVu_2V      \cr
I @>v_2>> V 
\endCD
$$
The projective cover $p\:PI \to I$ has as restriction
a surjective map $p'\:\rad PI \to \Omega S$. But $\rad PI$ is semisimple, thus 
$p'$ is a split epimorphism, thus we obtain a map $w\:\Omega S \to PI$
such that $pw = v_1.$ We consider the exact sequence induced from the sequence
$0 \to \Omega S \to PS \to S \to 0$ by the map $w$:
$$
\CD
0 @>>> \Omega S @>u_1>> PS @>e_1>> S @>>> 0 \cr
@.      @VwVV          @VVw'V      @| \cr
0 @>>> PI @>u_1'>> N @>e_1'>> S @>>> 0 
\endCD
$$
Here, $N$ is the pushout of the two maps $u_1$ and $w$. Since we know that $u_2u_1 = v_2v_1 = 
v_2pw$, there is a map $f\:N \to V$ such that $fu_1' = v_2p$ and $fw' = u_2$.
Since the map $\bmatrix v_2p & u_2\endbmatrix\:PI \oplus PS \to V$ is surjective, also $f$
is surjective.  

But recall that we assume that $\Ext^1(S,\Lambda) = 0,$ thus $\Ext^1(S,PI) = 0.$
This means that the lower exact sequence splits and therefore the socle of $N = PI \oplus S$
is a maximal submodule of $N$ (since $I$ is a local module, also $PI$ is a local module). 
Now $f$ maps the socle of $N$ into the socle of $V$, thus it maps a maximal submodule of $N$ 
into a simple submodule of $V$. This implies that the image of $f$ has length at most 2,
thus $f$ cannot be surjective. This contradiction shows that $Q$ has to be of length $2$, 
thus $Q = PS$ and therefore $PS$ is injective.
	\medskip
{\bf Lemma 3.}
 {\it If $S$ is a non-projective simple module with $\Ext^i(S,\Lambda) = 0$ for $1\le i \le d$,
then the modules $S_i = \Omega_iS$ with $0 \le i \le d$ are simple and not projective,
and the modules $P(S_i)$ are injective for $0 \le i < d.$}
	\medskip
The proof is by induction. If $d \ge 2$, we know by induction that the 
modules $S_i$ with $0 \le i \le d-1$ are simple and not projective,
and that the modules $P(S_i)$ are injective for $0 \le i < d-1.$
But $\Ext^1(\Omega_{d-1}S,\Lambda) \simeq \Ext^d(S,\Lambda) = 0,$ thus Lemma 2 asserts that
also $S_d$ is simple and not projective and that $P(S_{d-1})$ is 
injective.
	\medskip
{\bf Lemma 4.} {\it Let 
$S_0,S_1,\dots, S_b$ be simple modules with $S_i = \Omega_i(S_0)$ for
$1\le i \le b$. Assume that there is an integer $0 \le a < b$
such that the modules $S_i$ with $a\le i < b$ are pairwise non-isomorphic, whereas $S_b$ is isomorphic to $S_a$. In addition, we asssume that the modules $P(S_i)$
for $a \le i < b$ are injective. Then $S_a,\dots, S_{b-1}$ is the list of
all the simple modules and $\Lambda$ is self-injective.}
	\medskip
Proof: Let $\Cal S$ be the subcategory of all modules with composition factors of the
form $S_i$, where $a \le i < b.$ We claim that this subcategory is closed under 
projective covers and injective envelopes. Indeed, the projective cover of $S_i$
for $a \le i < b$ has the composition factors $S_i$ and $S_{i+1}$ (and $S_b = S_a$), thus
is in $\Cal S$.
Similarly, the injective envelope for $S_i$ with $a < i < b$ is $Q(S_i) = P(S_{i-1})$, thus it
has the composition factors $S_{i-1}$ and $S_i$, and $Q(S_a) = Q(S_b) = P(S_{b-1})$
has the composition factors $S_{b-1}$ and $S_a$.
Since we assume that $\Lambda$ is connected, we know that the only non-trivial subcategory closed under composition factors, extensions, projective covers and injective envelopes is the
module category itself. This shows that $S_a,\dots, S_{b-1}$ are all the simple modules.
Since the projective cover of any simple module is injective, $\Lambda$ is self-injective.
	\bigskip
{\bf Theorem 1.} {\it  Let $\Lambda$ be a connected left artinian ring with 
radical square zero. Assume that $\Lambda$ is not self-injective.
If $S$
is a non-projective simple module such that $\Ext^i(S,\Lambda) = 0$ for $1 \le i \le d$,
then the modules $S_i = \Omega_iS$ with $0 \le i \le d$ are pairwise 
non-isomorphic simple and non-projective modules and 
the modules $P(S_i)$ are injective for $0 \le i < d.$}
	\medskip
Proof. According to Lemma 3, the modules $S_i$ (with $0\le i \le d$) 
are simple and non-projective, 
and the modules $P(S_i)$ are injective for $0 \le i < d.$
If at least two of the modules $S_0, \dots, S_d$ are isomorphic, then 
Lemma 4 asserts that $\Lambda$ is self-injective, but this we have excluded.
	\bigskip
%=========================================================================
{\bf Theorem 2.} {\it Let $\Lambda$ be a connected left artinian ring with 
radical square zero and with $n$ simple modules. The following conditions are equivalent:}
\item{\rm (i)} {\it $\Lambda$ is self-injective, but not a simple ring.}
\item{\rm (ii)} {\it There exists a non-projective module $M$ with $\Ext^i(M,\Lambda) = 0$
 for $1\le i \le n+1$.}
\item{\rm (iii)} {\it There exists a non-projective simple module $S$ with $\Ext^i(S,\Lambda) = 0$
 for $1\le i \le n$.}
	\medskip
Proof. First, assume that $\Lambda$ is self-injective, but not simple. Since 
$\Lambda$ is not semisimple, there is a non-projective 
module $M$. Since $\Lambda$ is self-injective, $\Ext^i(M,\Lambda) = 0$ for all $i \ge 1.$
This shows the implication (i) $\implies$ (ii). 
The implication (ii) $\implies$ (iii)  follows from Lemma 1.
Finally, for the implication (iii) $\implies$ (ii) we use Theorem 1. 
Namely, if $\Lambda$ is not self-injective, then
Theorem 1 asserts that the simple modules $S_i = \Omega_iS$ with $0 \le i \le n$ are
pairwise non-isomorphic. However, these are $n+1$ simple modules, and we assume that the
number of isomorphism classes of simple modules is $n$. This completes the proof of Theorem 2.
	\bigskip
Note that the implication (ii) $\implies$ (i) in Theorem 2 asserts in particular that either 
$\Lambda$ is self-injective or else that any CM module
is projective, as shown by Chen [C]. 
Let us recall that a module $M$ is said to be a CM module
provided 
$\Ext^i(M,\Lambda) = 0$ and $\Ext^i(\Tr M,\Lambda) = 0$, for all $i \ge 1$ 
(here $\Tr$ denotes the transpose of the module); these modules are also called 
Gorenstein-projective modules, or totally reflexive modules, or modules of G-dimension equal to $0$. 
Note that in general there do exist modules $M$ with 
$\Ext^i(M,\Lambda) = 0$ for all $i \ge 1$ which are not CM modules, see [JS].
	\medskip
We also draw the attention to the generalized Nakayama conjecture 
formulated by Auslander-Reiten [AR]. It asserts
that for any artin algebra $\Lambda$ and any simple $\Lambda$-module $S$ 
there should exist an integer $i \ge 0$ such that 
$\Ext^i(S,\Lambda) \neq 0.$ It is known that this conjecture
holds true for algebras with radical square zero.  The implication (iii) $\implies$ (i) of Theorem 2
provides an effective bound: {\it If $n$ is the number of simple
$\Lambda$-modules, and $S$ is simple, then $\Ext^i(S,\Lambda) \neq 0$ for some $0 \le i \le n.$.}
Namely, in case $S$ is projective or $\Lambda$ is self-injective, 
then $\Ext^0(S,\Lambda) \neq 0.$ Now assume that $S$ is simple and not projective and that
$\Lambda$ is not self-injective. Then there must exist some integer 
$1\le i \le n$ with $\Ext^i(S,\Lambda) \neq 0,$ since otherwise the condition (iii) would be
satisfied and therefore condition (i).
	\bigskip
Theorem 1 may be interpreted as a statement concerning the $\Ext$-quiver 
of $\Lambda$. Recall that the {\it $\Ext$-quiver $\Gamma(R)$} of a left artinian ring $R$
has as vertices the (isomorphism classes of the)
simple $R$-modules, and if $S, T$ are simple $R$-modules, there is an arrow $T \to S$
provided $\Ext^1(T,S) \neq 0$, thus provided that there exists an 
indecomposable $R$-module $M$ of length 2 with socle $S$ and top $T$.
We may add to the arrow $\alpha\:T\to S$ the number $\l(\alpha) = ab$, where 
$a$ is the length of $\soc PT$ and $b$ is the length of
$QS/\soc$ (note that $b$ may be infinite).  
The properties of $\Gamma(R)$ which are relevant for this note 
are the following:
the vertex $S$ is a sink if and only if $S$ is projective; 
the vertex $S$ is a source if and only if $S$ is injective; finally, if $R$ is
a radical square zero ring and $S, T$ are simple $R$-modules then $PT = QS$ if and only if there
is an arrow $\alpha\:T \to S$ with $l(\alpha)=1$ and this is the only arrow starting
at $T$ and the only arrow ending in $S$.

Theorem 1 assert the following:
Let $\Lambda$ be a connected left artinian ring with 
radical square zero. Assume that $\Lambda$ is not self-injective.
Let $S$
be a non-projective simple module such that $\Ext^i(S,\Lambda) = 0$ for $1 \le i \le d$,
and let $S_i = \Omega_iS$ with $0 \le i \le d$. Then the local structure of
$\Gamma(\Lambda)$ is as follows:
$$
{\beginpicture
\setcoordinatesystem units <1.5cm,1cm>
\put{$S_0$} at 0 0  
\put{$S_1$} at 1 0
\put{$S_{d-1}$} at 2 0  
\put{$S_d$} at 3 0
\put{$\dots$} at 1.45 0
\arr{0.2 0}{0.8 0}
\arr{2.3 0}{2.8 0}

\arr{-.8 0.8}{-.2 .2}
\arr{-.8 0.4}{-.2 .05}
\arr{-.8 -.8}{-.2 -.2}

\arr{3.2 0.2}{3.8 .8}
\arr{3.2 0.05}{3.8 .4}
\arr{3.2 -.2}{3.8 -.8}
\multiput{$\vdots$} at -.7 -.1   3.7 -.1 /
\multiput{} at -1 0.8  4 -.8 /
\multiput{$1$} at 0.5 0.2  2.5 0.2 /
\endpicture}
$$ 
such that there is at least one arrow starting in $S_d$ (but may-be no arrow ending in
$S_0$). To be precise: the picture is supposed to 
show all the arrows starting or ending in the
vertices $S_0,\dots, S_d$ (and to assert that the vertices $S_0,\dots, S_d$
are pairwise different).
	\bigskip
Let us introduce the quivers $\Delta(n,t)$, where $n,t$ are positive integers. The quiver 
$\Delta(n,t)$ has $n$ vertices and also $n$ arrows, namely the
vertices labeled $0, 1,\dots, n-1$,
and arrows $i \to i\!+\!1$ for $0 \le i \le n-1$ (modulo $n$) (thus, we deal with an oriented
cycle); in addition, let $l(\alpha) = t$ for the arrow $\alpha\:n\!-\!1 \to 0$
and let $l(\beta) = 1$ for the remaining arrows $\beta$:
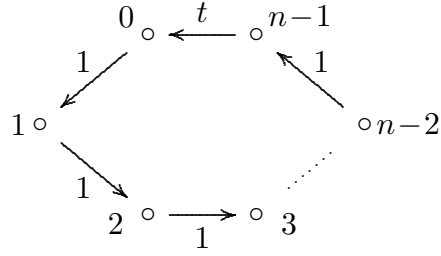
$$
{\beginpicture
\setcoordinatesystem units <1.2cm,1cm>
\multiput{$\circ$} at 0 0  1 1  1 -1  2 1  2 -1  3 0 /
\arr{0.8 0.8}{0.2 0.2}
\arr{0.2 -.2}{0.8 -.8}
%\arr{2.2 -.8}{2.8 -.2}
\arr{2.8 .2}{2.2 .8}
\arr{1.8 1}{1.2 1}
\arr{1.2 -1}{1.8 -1}
\setdots <1mm>
\plot 2.3 -.7  2.7 -.3 /
\setsolid
\put{$0$} at 0.8 1.2
\put{$1$} at -.2 0
\put{$2$} at 0.7 -1.1
\put{$3$} at 2.3 -1.1 
\put{$n\!-\!2$} at 3.4 0
\put{$n\!-\!1$} at 2.4 1.2
\multiput{$1$} at 0.4 0.7  0.4 -.7  1.5 -1.25  2.6 0.7   /
\put{$t$} at 1.5 1.25 
\endpicture}
$$
Note that the $\Ext$-quiver of a connected self-injective left artinian ring with
radical square zero and $n$ vertices is just $\Delta(n,1)$. Our further interest lies
in the cases  $t > 1.$ 
	\bigskip

{\bf Theorem 3.} {\it  Let $\Lambda$ be a connected left artinian ring with 
radical square zero and with $n$ simple modules.}

\item{(a)} {\it If there exists a non-projective simple modules $S$ with
  $\Ext^i(S,\Lambda) = 0$
 for $1\le i \le n-1,$
or if there exists a non-projective module $M$ with
  $\Ext^i(M,\Lambda) = 0$
 for $1\le i \le n$, 
 then $\Gamma(\Lambda)$ is of the form $\Delta(n,t)$ with $t > 1$.}
\item{(b)} {\it Conversely, if $\Gamma(\Lambda) = \Delta(n,t)$ and $t > 1$,
 then there exists a unique simple module $S$ with
  $\Ext^i(S,\Lambda) = 0$  for $1\le i \le n-1$, namely the module $S = S(0)$
(and it satisfies $\Ext^n(S,\Lambda) \neq 0)$.}
\item{(c)} {\it If $\Gamma(\Lambda) = \Delta(n,t)$ and $t > 1$,
 and if we assume in addition that $\Lambda$ is an artin algebra, then 
there exists a unique indecomposable module $M$ with
  $\Ext^i(M,\Lambda) = 0$
 for $1\le i \le n$, namely $M = \Tr D(S(0))$
  (and it satisfies $\Ext^{n+1}(M,\Lambda) \neq 0).$}
	\medskip
Here, for $\Lambda$ an artin algebra, 
$D$ denotes the $k$-duality, where $k$ is 
the center of $\Lambda$ (thus $D = \Hom_k(-,E)$, where $E$ is a minimal injective cogenerator
in the category of $k$-modules);
thus $D\Tr$ is the Auslander-Reiten translation and $\Tr D$ the reverse.
	\medskip
Proof of Theorem 3. Part (a) is a direct consequence of Theorem 1, using the
interpretation in terms of the $\Ext$-quiver as outlined above. Note that we must have
$t > 1$, since otherwise $\Lambda$ would be self-injective. 

(b) We assume that $\Gamma(\Lambda) = \Delta(n,t)$ with $t >1$.
For $0\le i < n$, let $S(i)$  be the simple module corresponding to the vertex $i$,
let $P(i)$ be its projective cover, $I(i)$ its injective envelope. We see from the quiver that
all the projective modules $P(i)$ with $0 \le i \le n\!-\!2$ are injective, thus
$\Ext^j(-,\Lambda) = \Ext^j(-,P(n\!-\!1))$ for all $j \ge 1.$
In addition, the quiver shows that 
$\Omega S(i) = S(i\!+\!1)$ for $0 \le i \le n\!-\!2$. Finally, we have $\Omega S(n\!-\!1) = S(0)^a$ for some positive integer
$a$ dividing $t$ and the injective envelope of $P(n\!-\!1)$ 
yields an exact sequence 
$$
 0 \to P(n\!-\!1) \to I(P(n\!-\!1)) \to S(n\!-\!1)^{t\!-\!1} \to 0 \tag{$*$}
$$
(namely, $I(P(n\!-\!1)) = I(\soc P(n\!-\!1)) = I(S(0)^a) = I(S(0))^a$ and 
$I(S(0))/\soc$ is the direct sum of $b$ copies of $S(n\!-\!1)$, where $ab= t$; 
thus the cokernel of the inclusion map $P(n\!-\!1) \to I(P(n\!-\!1))$ consists of
$t\!-\!1$ copies of $S(n\!-\!1)$).

Since $t>1$, the exact sequence $(*)$ shows that $\Ext^1(S(n\!-\!1),P(n\!-\!1)) \neq 0$.
It also implies that $\Ext^1(S(i),P(n\!-\!1)) = 0$ for $0 \le i \le n\!-\!2,$
and therefore that 
$$
\align
  \Ext^i(S(0),P(n\!-\!1)) &= \Ext^1(\Omega_{i-1}S(0),P(n\!-\!1)) \cr
     &= \Ext^1(S(i\!-\!1),P(n\!-\!1))\cr
   &=  0
\endalign
$$
for $1 \le i \le n\!-\!1$.

Since $\Omega_{n-i-1}S(i) = S(n\!-\!1)$ for $0\le i \le n\!-\!1$, we see that
$$
\align
  \Ext^{n-i}(S(i),P(n\!-\!1)) 
     &=   \Ext^{1}(\Omega_{n-i-1} S(i),P(n\!-\!1)) \cr
     &=   \Ext^{1}(S(n\!-\!1),P(n\!-\!1))\cr
     & \neq 0 
\endalign
$$
for $0 \le i \le n\!-\!1$.
Thus, on the one hand, we have $\Ext^n(S(0),\Lambda) \neq 0$, this
concludes the proof that $S(0)$ has the required properties. On the other hand,
we also see that $S = S(0)$ is the only simple module with 
$\Ext^i(S,\Lambda) = 0$ for $1\le i \le n\!-\!1$. This completes the proof of (b).

(c) Assume now in addition that $\Lambda$ is an artin algebra. As usual, we denote the
Auslander-Reiten translation $D\Tr$ by $\tau.$
Let $M$ be a non-projective
 indecomposable module with $\Ext^i(M,\Lambda) = 0$
 for $1\le i \le n$. The shape of $\Gamma(\Lambda)$
shows that $\Omega M = S^c$ for some simple module $S$ (and we have $c\ge 1$), also it shows
that no simple module is projective. 
Now $\Ext^i(S,\Lambda) = 0$ for $1 \le i < n$, thus according to (b) we must have
$S = S(0)$. It follows that
$PM$ has to be a direct sum of copies of $P(n\!-\!1)$, say of $d$ copies.
Thus a minimal projective presentation of
$M$ is of the form
$$
     P(0)^c \to P(n\!-\!1)^d \to M \to 0,
$$
and therefore a minimal injective copresentation of $\tau M$ is of the form
$$
  0 \to \tau M \to I(0)^c \to I(n\!-\!1)^d.
$$
In particular, $\soc \tau M = S(0)^c$ and 
$(\tau M)/\soc$ is a direct sum of copies of $S(n\!-\!1)$.

Assume that $\tau M \neq S(0)$, thus it has at least one composition factor of the form
$S(n\!-\!1)$ and therefore there exists a non-zero map $f\:P(n\!-\!1) \to \tau M.$
Since $\tau M$ is indecomposable and not injective, any map from an injective module to 
$\tau M$ maps into the socle of $\tau M$. But the image of $f$ is not contained in the socle of
$\tau M$, therefore $f$ cannot be factored through an injective module. It follows that
$$
 \Ext^1(M,P(n\!-\!1)) \simeq D\overline{\Hom}(P(n\!-\!1),\tau M) \neq 0,
$$
which contradicts the assumption that $\Ext^1(M,\Lambda) = 0.$ This shows that $\tau M = S(0)$
and therefore $M = \Tr D S(0).$

Of course, conversely we see that $M = \Tr D S(0)$ satisfies $\Ext^i(M,P(n\!-\!1)) = 0$
for $1\le i \le n,$ and $\Ext^{n+1}(M,P(n\!-\!1)) \neq 0$.
	\medskip
Remarks. (1) The module $M = \Tr D S(0)$ considered in (c) has length $t^2+t-1$, thus
the number $t$ (and therefore $\Delta(n,t)$) is determined by $M$.

(2) If $\Lambda$ is an artin algebra with $\Ext$-quiver $\Delta(n,t),$ the number $t$
has to be the square of an integer, say $t = m^2.$ A typical example of such an artin algebra
is the path algebra
of the following quiver
$$
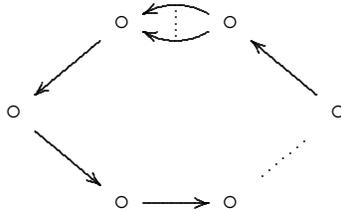
{\beginpicture
\setcoordinatesystem units <1.2cm,1cm>
\multiput{$\circ$} at 0 0  1 1  1 -1  2 1  2 -1  3 0 /
\arr{0.8 0.8}{0.2 0.2}
\arr{0.2 -.2}{0.8 -.8}
\arr{2.8 .2}{2.2 .8}
\arr{1.2 -1}{1.8 -1}
\setquadratic
\plot 1.8 1.1  1.5 1.2  1.2 1.1 /
\plot 1.8 .9  1.5 .8  1.2 .9 /
\arr{1.21 1.105}{1.2 1.1}
\arr{1.21 .895}{1.2 .9}

\setlinear
\setdots <.7mm>
\plot 1.5 0.85  1.5 1.15 /
\setdots <1mm>
\plot 2.3 -.7  2.7 -.3 /
\setsolid
\endpicture}
$$
with altogether $n+m-1$ arrows, modulo the ideal generated by all paths of length $2$.
Of course, if $\Lambda$ is a finite-dimensional $k$-algebra with radical square zero and 
$\Ext$-quiver $\Delta(n,m^2)$, and $k$ is an algebraically closed field, then
$\Lambda$ is Morita-equivalent to such an algebra. 

Also the following artin algebras with radical square zero and 
$\Ext$-quiver $\Delta(1,m^2)$ may be of interest: 
the factor rings of the polynomial ring $\Bbb Z[T_1,\dots,T_{m-1}]$
modulo the square of the ideal generated by
some prime number $p$ and the variables $T_1,\dots,T_{m-1}$.
	\bigskip\bigskip
{\bf References.}
	\medskip
\item{[AR]} M. Auslander, I. Reiten: On a generalized version of the Nakayama conjecture.
  Proc. Amer. Math. Soc. 52 (1975), 69--74.
\item{[ARS]} M. Auslander, I. Reiten, S. Smal\o{}: Representation Theory of Artin Algebras. 1995. 
   Cambridge University Press. 
\item{[C]} X.-W. Chen: Algebras with radical square zero are either self-injective or
 CM-free. Proc. Amer. Math. Soc. (to appear). 

\item{[JS]} D. A. Jorgensen, L. M. \c{S}ega: Independence of the total reflexivity conditions for
  modules. Algebras and Representation Theory 9 (2006), 217--226.
	\bigskip\bigskip

{\rmk
C. M. Ringel\par
Fakult\"at f\"ur Mathematik, Universit\"at Bielefeld, \par
POBox 100 131, D--33 501 Bielefeld, Germany, and \par
Department of Mathematics, Shanghai Jiao Tong University \par
Shanghai 200240, P. R. China. \par 
e-mail: {\ttk ringel\@math.uni-bielefeld.de} \par
\medskip

B.-L. Xiong\par
Department of Mathematics, Shanghai Jiao Tong University \par
Shanghai 200240, P. R. China.\par 
e-mail: {\ttk  xiongbaolin\@gmail.com} \par}
	\bye